\documentclass[a4paper,10pt]{article}
\pdfoutput=1
\usepackage[english]{babel}
\usepackage{amssymb}
\usepackage{amsmath}
\usepackage{amsthm}
\usepackage{graphicx}
\usepackage{fancyhdr}
\newfont{\bssten}{cmssbx10}
\newfont{\bssnine}{cmssbx10 scaled 900}
\newfont{\bssdoz}{cmssbx10 scaled 1200}

\usepackage{latexsym,amssymb,amsthm,amsxtra}
\usepackage{amsmath}
\usepackage{stmaryrd}
\usepackage{mathrsfs}
\usepackage{tikz} 
\usepackage{pgf} 
{\bf}{\it}
\newtheorem{theorem}{Theorem}

\newtheorem{lemma}{Lemma}
\newtheorem{remark}{Remark}
\newtheorem{proposition}{Proposition}

\newtheorem{ex}{Example}

\makeatletter
\DeclareRobustCommand{\cev}[1]{%
  \mathpalette\do@cev{#1}%
}
\newcommand{\do@cev}[2]{%
  \fix@cev{#1}{+}%
  \reflectbox{$\m@th#1\vec{\reflectbox{$\fix@cev{#1}{-}\m@th#1#2\fix@cev{#1}{+}$}}$}%
  \fix@cev{#1}{-}%
}
\newcommand{\fix@cev}[2]{%
  \ifx#1\displaystyle
    \mkern#23mu
  \else
    \ifx#1\textstyle
      \mkern#23mu
    \else
      \ifx#1\scriptstyle
        \mkern#22mu
      \else
        \mkern#22mu
      \fi
    \fi
  \fi
}

\makeatother

\def\epsilon{\varepsilon}

\usepackage{anysize}
\marginsize{3cm}{2.5cm}{2.5cm}{2.5cm} 
\newcommand{\be}{\begin{equation}}
\newcommand{\ee}{\end{equation}}

\newcommand{\norm}[1]{\left|#1\right|}
\newcommand\td[1]{\overline{#1}}

\def\R{{\mathbb R}}

\def\bV{{\mathbf V}}
\def\N{{\mathbb N}}

\def\T{{\mathcal{T}}}

\def\cv{\cev}

\def\bw{\mathbf w}
\def\mv{\mathbf v}
\def\mbB{{\mathbb B}}
\def\mbF{{\mathbb F}}

\def\maV{{\mathcal{V}}}

\def\maM{{\mathscr{M}}}

\def\maE{{\mathcal{E}}}
\def\maI{{\mathcal{I}}}
\def\maF{{\mathcal{F}}}

\def\I{{\mathbb{I}}}
\def\v{{\--}}
\def\pv{{\not\!\!\--}}

\def\T2a{{\tau_{2\vec \alpha^{+}}}}
\def\t2a{{t_{2\vec \alpha^{+}}}}
\newcommand\suite[1]{\left(#1\right)_{n\in\N}}

\def\\xi{{\cal{\xi}}}

\def\qed{\hfill $\blacksquare$}
\def\({{\Bigl(}}
\def\){{\Bigr)}}

\newcommand\pr[1]{{\mathbb P}\left[#1\right]}
\newcommand\esp[1]{{\mathbf E}\left[#1\right]}

\def\square{\ifmmode\sqr\else{$\sqr$}\fi}
\def\sqr{\vcenter{
         \hrule height.1mm
         \hbox{\vrule width.1mm height2.2mm\kern2.18mm\vrule width.1mm}
         \hrule height.1mm}}                  

{\end{list}}

\def\maM{\mathcal M}

\begin{document}
\title{A product form for the general stochastic matching model}
\author{P. Moyal, A. Bu$\check{\mbox{s}}$i\'c and J. Mairesse\\
\small{UTC/ Université de Lorraine, INRIA and CNRS/Universit\'e Pierre et Marie Curie}}
\maketitle

\begin{abstract}
\quad We consider a stochastic matching model with a general compatibility graph, as introduced in \cite{MaiMoy16}. 
We show that the natural necessary condition of stability of the system is also sufficient for the natural 
matching policy 'First Come, First Matched'  
(FCFM). For doing so, we derive the stationary distribution under a remarkable product form, by using an original dynamic reversibility property related to that of \cite{ABMW17} 
for the bipartite matching model. 
\end{abstract}

\section{Introduction}
\label{sec:intro}
Consider a general stochastic matching model (GM), as introduced in \cite{MaiMoy16}: 
items of various classes enter a system one by one, to be matched by couples and to leave the system immediately upon matching. Two items are compatible if and only if their classes 
are adjacent in a compatibility graph $G$ that is fixed beforehand. The classes of the entering items are drawn following a prescribed probability measure $\mu$ on $\maV$, the set 
of nodes of $G$. 
This model is a variant of the Bipartite Matching model (BM) introduced 
in \cite{CKW09}, see also \cite{AdWe}. In the BM, the compatibility graph is bipartite (say $\maV=\maV_1 \cup \maV_2$). Along 
the various natural applications of this model, the nodes of $\maV_1$ and $\maV_2$ represent respectively classes of {customers} and {servers}, 
{kidneys} and {patients}, blood {givers} and blood {receivers}, {houses} and {applicants}, and so on. 
The items are matched by couples of $\maV_1\times\maV_2$, and also {arrive} 
by couples of $\maV_1\times\maV_2$.  
The classes of the elements of the entering couples are random, and it is assumed in the aforementioned references that 
the class of the entering element of $\maV_1$ is always independent of that of the entering element of $\maV_2$.  
Several declinations of the BM can be found in a recent literature, along various applications: the Max-weight policy is shown asymptotically optimal for models with relaxations in \cite{BM14}. A model for taxi platforms is investigated in \cite{BC15,BC17}, in which the bipartite graph is bipartite-complete, but the nominal matchings between items are themselves random. A fluid model for an overloaded system representing an housing allocation platform is analyzed in \cite{TW08}, and \cite{BDPS11} addresses an organ transplant system under real-time constraint. An extension of GM models to the case where the matching structure is an hypergraph (i.e. a set of nodes with edges possibly gathering more 
than two nodes) was recently proposed in \cite{RM19}, addressing the form of the stability region in function of the geometry of the considered hypergraph. 

An important generalization of the BM is the so-called {Extended Bipartite Matching} model (EBM) introduced in \cite{BGMa12}, where this independent assumption is relaxed. 
Possible entering couples are element of a bipartite {arrival graph} on the bipartition $\maV_1\cup\maV_2$. The construction of perfect infinite matchings for EBM models, i.e. of random graphs gathering the entered nodes and whose edges 
represent the executed matches, was recently presented in \cite{MBM18}. 
Importantly, notice that the GM can in fact be seen as a particular case of EBM, taking the bipartite double cover of $G$ as compatibility graph, and duplicating arrivals 
with a copy of an item of the same class. 

The main question raised in \cite{MaiMoy16} is the shape of the {stability region} of the model, that is, the set of probability measures on $\maV$ rendering the natural Markov chain of corresponding system positive recurrent. Partly relying on the aforementioned connection between GM and EBM, and the results of \cite{BGMa12}, \cite{MaiMoy16} show that the stability region is always included in the set of measures satisfying the natural condition (\ref{eq:Ncond}) below. The form of the stability region is then heavily dependent on the 
structural properties of the compatibility graph, and on the {matching policy}, i.e. the rule of choice of a match for an entering item whenever several matches are possible. A matching policy is then said to have a {maximal} stability region for $G$ if the system is stable for any measure satisfying (\ref{eq:Ncond}). It is shown in \cite{MaiMoy16} that a GM on a bipartite $G$ is never stable, that a designated class of graphs 
(the complete $k$-partite graphs for $k\ge 3$, see below) makes the stability region maximal for all matching policies, and that the policy 'Match the Longest' always has a maximal stability region for a non-bipartite $G$. 
Applying fluid unstability arguments to a continuous-time version of the GM, 
\cite{MoyPer17} shows that, aside for a very particular class of graphs, whenever $G$ is not complete $k$-partite there {always} exists a policy of the strict priority type that does not have a maximal stability region, and that the 'Uniform' random policy (natural in the case where no information is available to the entering items on the state of the system) never has a maximal stability region, 
thereby providing a partial converse of the result in \cite{MaiMoy16}. Related models are studied in 
\cite{GurWa} and \cite{NS17}, proposing optimization schemes for models of various matching structures (general graphs and hypergraphs), and the executions of matchings (which are possibly retarded, i.e. compatible items may not to be matched right away in order to wait for a more profitable future match) are associated to a cost or a reward. 
 
In this paper we are concerned with the stability region of the GM under the `First Come, First Matched' (FCFM) policy, consisting in always performing the match of the entering item with the oldest compatible item in the system. 
Compared to the aforementioned stability studies, this matching policy raises technical problems, mainly due to the intricate form of the state space of its natural Markov representation. In the history of study of the BM,  
the corresponding First Come, First Served' (FCFS) policy was the first one considered in the seminal papers \cite{CKW09,AdWe}, which show the existence of a stationary matching whenever the measure $\mu$ 
satisfies a natural resource pooling condition analog to the one in (\ref{eq:Ncond}) below. Further, \cite{ABMW17} recently showed that the stationary distribution of the Markov chain of the system can be obtained in a remarkable product form, which is obtained using an original dynamic reversibility argument. However, these results cannot be directly applied to the present context, for the GM is {\em not} a particular case of a BM, but of an EBM, for which the latter reversibility argument is unlikely to hold in general. (Observe that the maximality of FCFS for the EBM is conjectured, but left as on open problem in \cite{BGMa12}.) We show hereafter that we can in fact construct a reversibility scheme that is related to the one proposed in \cite{ABMW17} for the BM, to show that the stability region of FCFM is indeed maximal for the GM, and that the stationary state of the Markov representation also satisfies a product form. 

The paper is organized as follows. In Section \ref{sec:model} we introduce and formalize our model. 
In Section \ref{sec:FCFM} we develop our reversibility scheme for the FCFM system, leading to our main result, Theorem \ref{thm:main}, which establishes the existence of a stationary probability under a product form for the natural Markov representation of the system, under the natural condition that $\mu$ is an element of the set defined by (\ref{eq:Ncond}). Under the same condition, the construction of perfect infinite FCFM-matchings is presented in Section \ref{FCFMreverse}.

\section{The model}
\label{sec:model}

\subsection{General notation}
\label{subsec:notation} 
Denote by $\R$ the set of real numbers, by $\N$ the set of non-negative integers and by $\N_+$, the subset of positive integers. For any two integers $m$ and $n$, denote by $\llbracket m,n \rrbracket=[m,n] \cap \N$. 
For any $n\in\N_+$, let $\mathfrak S_n$ be the group of permutations of $\llbracket 1,n \rrbracket$. 
Let $A^*$ (respectively, $A^{\N}$) be the set of finite (resp., infinite) words over the alphabet $A$. 
Denote by $\emptyset$, the empty word of $A^*$. 
For any word $w \in A^*$ and any subset $B$ of $A$, we 
let $\norm{w}_B$ be the number of occurrences of elements of $B$ in $w$. 
For any  letter $a\in A$, we denote 
$\norm{w}_a:=\norm{w}_{\{a\}}$, and for any finite word $w$ we let $\norm{w}=\sum_{a\in A} \norm{w}_a$ be the
{\em length} of $w$. For a word $w \in A^*$ of length $\norm{w}=q$, we write $w=w_1w_2...w_q$, i.e. $w_{i}$ is the $i$-th letter 
of the word $w$. In the proofs below, we understand the word $w_1...w_k$ as $\emptyset$ whenever $k=0$. 
Also, for any $w \in A^*$ and any $i\in \llbracket 1,\norm{w} \rrbracket$, we denote by $w_{[i]}$, the word of length $\norm{w}-1$ obtained from $w$ by deleting 
its $i$-th letter. 

Consider a simple undirected graph $G=(\maV,\maE)$, where $\maV$ denotes the set of nodes, and $\maE \subset \maV^2\setminus \Delta$ is the set of edges, where $\Delta=\left\{(i,i):\,i\in \maV\right\}$, i.e. we do not allow self-loops. 
We use the
notation $u \v v$ for $(u,v) \in \maE$
and $u \pv v$ for $(u,v) \not\in \maE$. 
For $U \subset \maV$, we define $U^c = \maV \setminus U$ and
\[
\maE(U) = \{v \in \maV\,:\, \exists u \in U, \ u
\-- v\}\:.
\]
An {\em independent set} of $G$ is a non-empty subset $\maI \subset \maV$ 
which does not include any pair of neighbors, {\em i.e.} $\bigl(\forall i, j \in \maV, \, i \v j \Rightarrow i\not\in \maI \mbox{ or }j\not\in \maI\bigr)$. 
Let $\I(G)$ be the set of independent sets of $G$. An independent set $\maI$ is said to be {\em maximal} if $\maI \cup \{j\} \not\in \I(G)$ for any $j \not\in \maI$. 

Recall that a connected graph $G=(\maV,\maE)$ is said {\em complete $p$-partite}, 
$p\ge 2$, if there exists a partition of $\maV$ into maximal independent sets $\maI_1,\dots, \maI_p$, such that
 \[
 \forall i\neq j,\, \forall u \in \maI_i,\, \forall v \in \maI_j,\,\,u \v v \:.
 \]
 Notice that such graphs are referred to as {\em separable} in the references \cite{MaiMoy16} and \cite{MoyPer17}. 

\medskip

For $w=w_1w_2\,...\,w_q \in \maV^*$, we denote by $\cv w$ the reversed version of $w$, i.e. $\cv w=w_qw_{q-1}...w_2w_1.$ 
Let $\td\maV$ be an independent copy of the set $\maV$, i.e. a set of cardinality $|\maV|$, disjoint of $\maV$ and containing 
copies of the elements of $\maV$. We denote by $\td a$, the copy of any element $a$ of $\maV$, and also say that $\td a$ is the {\em counterpart} of $a$, and vice-versa. 
For any $\td a \in \td \maV$, let us denote $\td{\td a}=a$. 
Let $\bV:=\maV \cup \td\maV.$ For any word $\bw \in \mathbf V^*$, denote by $\maV(\bw)$ (respectively, $\td\maV(\bw)$) the set of letters of $\maV$ (resp., $\td\maV$) that are present in $\bw$: 
\begin{equation*}
\maV(\bw) =\Bigl\{a \in \maV\,:\,\norm \bw_a >0\Bigl\};\quad \quad\td\maV(\bw) =\Bigl\{\td a \in \td\maV\,:\,\norm \bw_{\td a} >0\Bigl\}.
\end{equation*}
For any $\bw \in \bV^*$, the {\em restriction} of $\bw$ to $\maV$ (respectively, to $\td\maV $) is the word 
$\bw|_{\maV} \in \maV^*$ (resp., $\bw|_{\td\maV} \in \td\maV^*$) of size $\norm{\bw}_\maV$ (resp. of size $\norm{\bw}_{\td\maV}$), 
obtained by keeping only the letters belonging to $\maV$ (resp. to $\td\maV$) in $\bw$, in the same order. 
The {\em dual} $\td \bw$ of the word $\bw=\bw_1...\bw_q \in \bV^*$ is the word obtained 
by exchanging the letters of $\bw$ with their counterpart, i.e. $\td \bw=\td{\bw_1}\,...\,\td{\bw_q}.$ 

\begin{ex}\rm
Take for instance  
$\bw=a\,b\,\td a\,c\,\td b\,\td c\,\td b\,d\,a$. Then we obtain 
$\maV(\bw) =\left\{a,b,c,d\right\},\,\,\td\maV(\bw) =\left\{\td a, \td b,\td c\right\}$,\,\, $\bw|_{\maV} =a\,b\,c\,d\,a,\,\,\bw|_{\td\maV}=\td a\,\td b\,\td c\,\td b,\,\,
\td \bw =\td a\,\td b\,a\,\td c\,b\,c\,b\,\td d\,\td a,\,\,\cv \bw = a\,d\,\td b\,\td c\,\td b\,c\,\td a\,b\,a.$
\end{ex}

\subsection{Formal definition of the model}
\label{subsec:model}
We consider a {\em general stochastic matching model}, as was defined in \cite{MaiMoy16}: items enter one by one a system, and each of them belongs to a 
determinate class. The set of classes is denoted by $\maV$, and identified with $\llbracket 1,|\maV|\rrbracket$. We fix a connected simple graph $G=(\maV,\maE)$ having set of nodes $\maV$, termed {\em compatibility graph}. 
Upon arrival, any incoming item of class, say, $i \in \maV$ is either matched with an item present in the buffer, of a class $j$ such that 
$i \v j$, if any, or if no such item is available, it is stored in the buffer to wait for its match. 
Whenever several matches are possible for an incoming item $i$, it is the role of the matching policy to decide the match of the latter item without ambiguity. 
Throughout, we assume that the matching policy is 'First Come, First Matched' (\textsc{fcfm}), that is, the match of $i$ is the oldest among all stored items of neighboring classes of $i$. Each matched pair departs the system right away. 

We assume that the successive classes of entering items are random.  We fix an underlying probability space $(\Omega,\mathcal F,\mathbb P)$, on which all random variables (r.v.'s, for short) are defined. For any $n \in \N$, let $V_n \in \maV$ denote the class of the $n$-th incoming item. Throughout this work, we assume that the sequence $\suite{V_{n}}$ is iid, from the distribution $\mu$ on $\mathcal V$. Without loss of generality, we assume that $\mu$ has full support $\maV$ 
(we write $\mu \in \maM(\maV)$). Altogether, according to the terminology in \cite{MaiMoy16}, we thus consider the general matching (GM) model associated to 
$(G,\mu,\textsc{fcfm})$.

\subsection{Markov representation}
\label{subsec:Markov}
Fix the compatibility graph $G=(\maV,\maE)$ until the end of this section. 
Fix an integer $n_0 \ge 1$, and a realization $v_1,...,v_{n_0}$ of $V_1,...,V_{n_0}$. 
Define the word $\mv\in \maV^*$ by $\mv:=v_1...v_{n_0}$. 
Then, there exists a unique {\sc fcfm} {\em matching} of the word $\mv$, that is, a graph having set of nodes 
$\left\{v_1,...,v_{n_0}\right\}$ and whose edges represent the matches performed in the system until time $n_0$, if the successive arrivals are given by $\mv$ 
and the matching policy is {\sc fcfm}. 
This matching is denoted by $M^{\textsc{fcfm}}(\mv)$. 
The state of the system is defined as the word $W^{\textsc{fcfm}}(\mv)\in \maV^*$, whose letters are the classes of the unmatched items at time $n_0$, 
i.e. the isolated vertices in the matching $M^{\textsc{fcfm}}(\mv)$, in their order of arrivals. The word $W^{\textsc{fcfm}}(\mv)\in \maV^*$ is called 
{\em queue detail} at time $n_0$. 
Then any admissible queue detail belongs to the set 
\begin{equation}
\mathbb W = \Bigl\{ w\in \maV^*\; : \; \forall  (i,j) \in \maE, \; |w|_i|w|_j=0  \Bigr\}.\label{eq-ncss}
     \end{equation} 

Fix a (possibly random) word $Y \in \mathbb W$. Denote for all $n\ge 0$ by $W^{\{Y\}}_n$ the buffer content at time $n$ 
(i.e. just before the arrival of item $n$) if the buffer content at time 0 was set to $w$, in other words 
\[W^{\{Y\}}_n= W^\phi\left(YV_0...V_{n-1}\right).\]
Then the buffer-content sequence is stochastic recursive, since we clearly have that 
\[\left\{\begin{array}{ll}
W^{\{Y\}}_0 &= Y;\\
W^{\{Y\}}_{n+1} &=W^{\{Y\}}_n \odot_{\textsc{fcfm}} (V_n),\,n\in\N,
\end{array}\right.\]
where for all $w\in\mathbb W$ and $v\in\maV$, 
$$w \odot_{\textsc{fcfm}} (v) =
\left \{
\begin{array}{ll}
wv & \textrm{if } \; |w|_{\maE(v)} = 0;\\
w_{\left [\Phi(w,v)\right]} & \textrm{else, where }\Phi(w,v) = \min \{k\in \llbracket 1,|w| \rrbracket,\,:\, w_k\in\maE(v)\}.
\end{array}
\right .
$$
Consequently, if we assume that the sequence $\suite{V_n}$ is independent of $Y$, the queue detail $\suite{W^{\{Y\}}_n}$ is a $\mathbb W$-valued 
$\mathcal F_n$-Markov chain, where $\mathcal F_0 = \sigma(Y)$ and $\mathcal F_n=\sigma\left(Y,V_0,...,V_{n-1}\right)$ for all $n\ge 1$. 
The sequence $\suite{W_n}$ is termed {\em natural chain} of the system.

\subsection{Stability of the matching model}
\label{subsec:stab}
Fix a connected graph $G$. The chain $\suite{W_n}$ is clearly irreducible, as all states of $\mathbb W$ lead to $\emptyset$.  
In line with \cite{MaiMoy16} we define the {\em stability region} of the GM model $(G,\mu,\textsc{fcfm})$ by 
\begin{equation}
\label{eq:defstab}
\textsc{stab}(G,\textsc{fcfm}) := \left\{\mu \in \maM\left(\maV\right)\,:\,\suite{W_n} \mbox{ is positive recurrent}\right\},
\end{equation}
which is clearly independent of the initial state $Y$ in view of the above observation.  
Let us also define the set of measures 
\begin{equation}
\label{eq:Ncond}
\textsc{Ncond}(G): \left\{\mu \in \maM\left(\maV\right)\,:\,\mbox{for any }\maI \in \I(G), \mu\left(\maI\right) < \mu\left(\maE\left(\maI\right)\right)\right\}
\end{equation}
which, from Theorem 1 in \cite{MaiMoy16}, is non-empty if and only if $G$ is non-bipartite. 
We know from Proposition 2 in [{\em ibid.}] that for any policy $\Phi$ in the set of so-called {\em admissible} matching policies (that includes {\sc fcfm}), 
the stability region of the GM model associated to $(G,\mu,\Phi)$, defined similarly to (\ref{eq:defstab}), is included in \textsc{Ncond}$(G)$. 
An admissible policy $\Phi$ is said to be {\em maximal} if these two sets coincide. Theorem 2 in [{\em ibid.}] establishes the maximality 
of the matching policy 'Match the Longest' for any non-bipartite graph, however priority policies and the uniform policy are in general not maximal (respectively, Theorem 3 and Proposition 7 in \cite{MoyPer17}). Last, it is stated in Theorem 2 in \cite{MaiMoy16} that any GM model on a graph $G$ that is $p$-partite complete for 
$p\ge 3$ has a stability region that coincides with \textsc{Ncond}($G$) whatever the admissible matching policy; in other word {\em any} admissible policy is maximal in this case. 

We show in this paper that the policy {\sc fcfm} is maximal, and characterize the steady state of the system under the condition that $\mu \in \textsc{Ncond}(G)$.

\section{Product form}
\label{sec:FCFM}
In this section we establish our main result, Theorem \ref{thm:main}, which establishes the maximality of the policy First Come, First Matched, by  
constructing explicitly the stationary distribution of the natural chain on $\mathbb W$. Interestingly enough, this probability distribution has a remarkable product form, detailed 
in (\ref{eq:PiW}). 
\begin{theorem}
\label{thm:main}
Let $G=(\maV,\maE)$ be a non-bipartite graph. Then the sets $\textsc{stab}(G,\textsc{fcfm})$ and $\textsc{Ncond}(G)$, defined respectively by (\ref{eq:defstab}) and (\ref{eq:Ncond}), coincide. 
In other words the GM model 
$(G,\mu,\textsc{fcfm})$ is stable if and only if $\mu$ belongs to the set \textsc{Ncond}($G$). In that case, the following is the only stationary probability of the natural chain $\suite{W_n}$:
\begin{equation}
\label{eq:PiW}
\Pi_W\left(w\right)=\alpha\prod\limits_{\ell =1}^q {\mu(w_\ell) \over \mu\Bigl(\maE\left(\left\{w_1,...,w_\ell\right\}\right)\Bigl)},\,
\mbox{ for any }w=w_1...w_q \in \mathbb W,
\end{equation}
where 
\begin{equation}
\label{eq:alpha}
\alpha =  \left\{1 +\sum_{\maI\in\mathbb I(G)} \sum_{\sigma\in\mathfrak S_{|\maI|}} \prod_{j=1}^{|\maI|} 
{\mu(i_{\sigma(j)}) \over \mu\left(\maE\left(\{i_{\sigma(1)},...,i_{\sigma(j)}\}\right)\right)-\mu\left(\left\{i_{\sigma(1)},...,i_{\sigma(j)}\right\}\right)}
\right\}^{-1}.
\end{equation}
\end{theorem}
The remainder of this Section is devoted to the proof of Theorem \ref{thm:main}, which, as will be demonstrated below, is based on a subtle reversibility 
scheme that is related to the proof of reversibility for the BM model in \cite{ABMW17}. Let us recall however that the GM model is not a particular case of BM model, 
so the proof below presents many specificities with respect to \cite{ABMW17}.

\subsection{Auxiliary Markov representations}
\label{subsec:Markov} 
We now introduce two auxiliary Markov representations of the system: the $\mathbf V^*$-valued Backwards and Forwards detailed chains, similar in construction to the backwards and forwards 'pair by pair detailed FCFS matching processes', introduced in subsection 5.1 of \cite{ABMW17}. 
\paragraph{Backwards detailed chain} 
We define the $\bV^*$-valued {backwards detailed process} $\suite{B_n}$ as follows: 
$B_0=\emptyset$ and for any $n\ge 1$, 
\begin{itemize}
\item if $W_n=\emptyset$ (i.e. all the items arrived up to time $n$ are matched at time $n$), then we set $B_n=\emptyset$;
\item if not, we let $i(n)\le n$ be the index of the oldest item in line.  
Then, the word $B_n$ is of length $n-i(n)+1$, and for any $\ell \in \llbracket 1,n-i(n)+1 \rrbracket$, we set 
\[B_n(\ell)=\left\{\begin{array}{ll}
V_{i(n)+\ell-1} \,\, &\mbox{if $V_{i(n)+\ell-1}$ has not been matched up to time $n$};\\
\td{V_{k}}\,\, &\mbox{if $V_{i(n)+\ell-1}$ is matched at or before time $n$, with item $V_k$}\\
               &\mbox{(where $k \le n$)}.\\
\end{array}\right.\] 
\end{itemize}
In other words, assuming that the initial system is empty, 
the word $B_n$ gathers the class indexes of all unmatched items entered up to $n$, and the copies of the class  
indexes of the items matched after the arrival of the oldest unmatched item at $n$, at the place of the class index of the item they have 
been matched to. 
Observe that we necessarily have that $B_n(1)=V_{i(n)}\in \maV$. Moreover, 
the word $B_n$ necessarily contains all the letters of $W_n$. More precisely, we have 
\begin{equation}
\label{eq:WB}
W^{\{\emptyset\}}_n=B_n|_{\maV},\,\,\,n\ge 0.
\end{equation}
It is easily seen that $\suite{B_n}$ also is a $\maF_n$-Markov chain: for any $n\ge 0$, the value of $B_{n+1}$ can be deduced from that of 
$B_n$ and the class $V_{n+1}$ of the item entered at time $n+1$.

\paragraph{Forward detailed chain} 
The $\bV^*$-valued {forward detailed process} $\suite{F_n}$ is defined as follows: 
$F_0=\emptyset$ and for any $n\ge 1$, 
\begin{itemize}
\item if $W_n=\emptyset$, then we also set $F_n=\emptyset$;
\item if not, we let $\mathscr U_n$ be the set of items entered up to $n$ and not yet matched at $n$ (which is non empty since $W_n \ne \emptyset$), and set 
      \[j(n)=\sup\left\{ m > n\,:\, \mbox{$V_m$ is matched with an element of $\mathscr U_n$}\right\}.\]
      Notice that $j(n)$ is possibly infinite. Then, $F_n$ is the word of $\bV^*$ of size $j(n)-n$ (respectively of $\bV^{\N}$ if $j(n)=+\infty$), such that for any $\ell \in \llbracket 1,j(n)-n \rrbracket$ 
      (resp., $\ell \in \N_+$), 
      \[F_n(\ell)=\left\{\begin{array}{ll}
                         V_{n+\ell} \,\, &\mbox{if $V_{n+\ell}$ is not matched with an item arrived up to $n$};\\
                        \td{V_{k}}\,\, &\mbox{if $V_{n+\ell}$ is matched with item $V_k$, where $k \le n$}.
                        \end{array}\right.\] 
\end{itemize}
In other words, assuming that the initial system is empty, the word $F_n$ contains the copies of all the class indexes of the items entered up to time $n$ and matched after $n$, 
together with the class indexes of all unmatched items entered before the last item matched with an item entered up to $n$, if any.  
Observe that if $F_n \in \bV^*$ is finite, then $F_n\left(j(n)-n\right) \in \td\maV$ since by definition, the item $V_{j(n)}$ is matched with some $V_k$ for $k\le n$, and therefore $F_n\left(j(n)-n\right)=\td V_k$. 
It is also clear that $\suite{F_n}$ is a $\maF_{n}$-Markov chain, as for any $n\ge 0$, the value of $F_{n+1}$ depends solely on $F_n$ 
and the class index $V_{j(n)+1}$ of the item entered at time $n+j(n)+1$. 


\begin{ex}
\label{ex:trajBWF}
\emph{Consider the compatibility graph of Figure \ref{fig:example1}, addressed in Section 5 of \cite{MaiMoy16} 
(this is the smallest graph that is neither bipartite nor complete $p$-partite). An arrival scenario together with successive values of the 
natural chain, the backwards and the forwards chain are represented in Figure \ref{fig:example2}.} 
\begin{figure}[h!]
\begin{center}
\begin{tikzpicture}[scale=1]
\draw[-] (2,3) -- (2,2);
\draw[-] (2,2) -- (1,1);
\draw[-] (2,2) -- (3,1);
\draw[-] (1,1) -- (3,1);
\fill (2,3) circle (2pt) node[right] {\small{1}} ;
\fill (2,2) circle (2pt) node[right] {\small{2}} ;
\fill (1,1) circle (2pt) node[below] {\small{3}} ;
\fill (3,1) circle (2pt) node[below] {\small{4}} ;
\end{tikzpicture}
\vspace*{-0.3cm}
\caption[smallcaption]{Compatibility graph of Example \ref{ex:trajBWF}.}
\label{fig:example1}
\end{center}
\end{figure}

\begin{figure}[h!]
\begin{center}
\begin{tikzpicture}[scale=0.7]
\draw[-] (-2,13) -- (19,13);
\draw[-] (-2,13) -- (-2,-2);
\draw[-] (-2,-2) -- (19,-2);
\draw[-] (19,-2) -- (19,13);
\draw[-,very thick] (-0.7,11.3) -- (-0.7,11.7);
\draw[-] (-1,11.5) -- (11,11.5);
\fill (0,11.5) circle (2pt) node[below] {\small{1}};
\fill (1,11.5) circle (2pt) node[below] {\small{3}};
\fill (2,11.5) circle (2pt) node[below] {\small{4}};
\fill (3,11.5) circle (2pt) node[below] {\small{$2$}};
\fill (4,11.5) circle (2pt) node[below] {\small{3}};
\fill (5,11.5) circle (2pt) node[below] {\small{1}};
\fill (6,11.5) circle (2pt) node[below] {\small{3}};
\fill (7,11.5) circle (2pt) node[below] {\small{2}};
\fill (8,11.5) circle (2pt) node[below] {\small{2}};
\fill (9,11.5) circle (2pt) node[below] {\small{1}};
\fill (10,11.5) circle (2pt) node[below] {\small{4}};
\fill (11,11.5) node[right]{\small{$W_0 =\emptyset,\quad B_0=\emptyset,\quad F_0=\emptyset$}};
 %
\draw[-] (-1,10) -- (11,10);
\fill (0,10) circle (2pt) node[below] {\small{1}};
\draw[->, thin] (0,10) .. controls +(up:0.5cm)  .. (3,10);
\draw[-,very thick] (0.3,9.8) -- (0.3,10.2);
\fill (1,10) circle (2pt) node[below] {\small{3}};
\fill (2,10) circle (2pt) node[below] {\small{4}};
\fill (3,10) circle (2pt) node[below] {\small{$\not 2 \, \bar 1$}};
\fill (4,10) circle (2pt) node[below] {\small{3}};
\fill (5,10) circle (2pt) node[below] {\small{1}};
\fill (6,10) circle (2pt) node[below] {\small{3}};
\fill (7,10) circle (2pt) node[below] {\small{2}};
\fill (8,10) circle (2pt) node[below] {\small{2}};
\fill (9,10) circle (2pt) node[below] {\small{1}};
\fill (10,10) circle (2pt) node[below] {\small{4}};
\fill (11,10) node[right]{\small{$W_1 =1,\quad B_1=1,\quad F_1=34\bar 1$}};
 %
\draw[-] (-1,8.5) -- (11,8.5);
\fill (0,8.5) circle (2pt) node[below] {\small{1}};
\draw[->, thin] (0,8.5) .. controls +(up:0.5cm)  .. (3,8.5);
\fill (1,8.5) circle (2pt) node[below] {\small{3}};
\draw[->, thin] (1,8.5) .. controls +(up:0.5cm)  .. (2,8.5);
\draw[-,very thick] (1.3,8.3) -- (1.3,8.7);
\fill (2,8.5) circle (2pt) node[below] {\small{$\not 4\,\bar 3$}};
\fill (3,8.5) circle (2pt) node[below] {\small{$\not 2 \, \bar 1$}};
\fill (4,8.5) circle (2pt) node[below] {\small{3}};
\fill (5,8.5) circle (2pt) node[below] {\small{1}};
\fill (6,8.5) circle (2pt) node[below] {\small{3}};
\fill (7,8.5) circle (2pt) node[below] {\small{2}};
\fill (8,8.5) circle (2pt) node[below] {\small{2}};
\fill (9,8.5) circle (2pt) node[below] {\small{1}};
\fill (10,8.5) circle (2pt) node[below] {\small{4}};
\fill (11,8.5) node[right]{\small{$W_2 =13,\quad B_2=13,\quad F_2=\bar 3\bar 1$}};
 %
\draw[-] (-1,7) -- (11,7);
\fill (0,7) circle (2pt) node[below] {\small{1}};
\draw[->, thin] (0,7) .. controls +(up:0.5cm)  .. (3,7);
\fill (1,7) circle (2pt) node[below] {\small{$\not 3\,\bar 4$}};
\draw[->, thin] (1,7) .. controls +(up:0.5cm)  .. (2,7);
\fill (2,7) circle (2pt) node[below] {\small{$\not 4\,\bar 3$}};
\draw[-,very thick] (2.3,6.8) -- (2.3,7.2);
\fill (3,7) circle (2pt) node[below] {\small{$\not 2 \, \bar 1$}};
\fill (4,7) circle (2pt) node[below] {\small{3}};
\fill (5,7) circle (2pt) node[below] {\small{1}};
\fill (6,7) circle (2pt) node[below] {\small{3}};
\fill (7,7) circle (2pt) node[below] {\small{2}};
\fill (8,7) circle (2pt) node[below] {\small{2}};
\fill (9,7) circle (2pt) node[below] {\small{1}};
\fill (10,7) circle (2pt) node[below] {\small{4}};
\fill (11,7) node[right]{\small{$W_3 =1,\quad B_3=1\bar 4\bar 3,\quad F_3=\bar 1$}};
 %
\draw[-] (-1,5.5) -- (11,5.5);
\fill (0,5.5) circle (2pt) node[below] {\small{1}};
\draw[->, thin] (0,5.5) .. controls +(up:0.5cm)  .. (3,5.5);
\fill (1,5.5) circle (2pt) node[below] {\small{3}};
\draw[->, thin] (1,5.5) .. controls +(up:0.5cm)  .. (2,5.5);
\fill (2,5.5) circle (2pt) node[below] {\small{4}};
\fill (3,5.5) circle (2pt) node[below] {\small{2}};
\draw[-,very thick] (3.3,5.3) -- (3.3,5.7);
\fill (4,5.5) circle (2pt) node[below] {\small{3}};
\fill (5,5.5) circle (2pt) node[below] {\small{1}};
\fill (6,5.5) circle (2pt) node[below] {\small{3}};
\fill (7,5.5) circle (2pt) node[below] {\small{2}};
\fill (8,5.5) circle (2pt) node[below] {\small{2}};
\fill (9,5.5) circle (2pt) node[below] {\small{1}};
\fill (10,5.5) circle (2pt) node[below] {\small{4}};
\fill (11,5.5) node[right]{\small{$W_4=\emptyset,\quad B_4=\emptyset,\quad F_4=\emptyset$}};
 %
\draw[-] (-1,4) -- (11,4);
\fill (0,4) circle (2pt) node[below] {\small{1}};
\draw[->, thin] (0,4) .. controls +(up:0.5cm)  .. (3,4);
\fill (1,4) circle (2pt) node[below] {\small{3}};
\draw[->, thin] (1,4) .. controls +(up:0.5cm)  .. (2,4);
\fill (2,4) circle (2pt) node[below] {\small{4}};
\fill (3,4) circle (2pt) node[below] {\small{2}};
\fill (4,4) circle (2pt) node[below] {\small{3}};
\draw[-,very thick] (4.3,3.8) -- (4.3,4.2);
\draw[->, thin] (4,4) .. controls +(up:0.5cm)  .. (7,4);
\fill (5,4) circle (2pt) node[below] {\small{1}};
\fill (6,4) circle (2pt) node[below] {\small{3}};
\fill (7,4) circle (2pt) node[below] {\small{$\not 2\,\bar 3$}};
\fill (8,4) circle (2pt) node[below] {\small{2}};
\fill (9,4) circle (2pt) node[below] {\small{1}};
\fill (10,4) circle (2pt) node[below] {\small{4}};
\fill (11,4) node[right]{\small{$W_5 =3,\quad B_5=3,\quad F_5=13\bar 3$}};
 %
\draw[-] (-1,2.5) -- (11,2.5);
\fill (0,2.5) circle (2pt) node[below] {\small{1}};
\draw[->, thin] (0,2.5) .. controls +(up:0.5cm)  .. (3,2.5);
\fill (1,2.5) circle (2pt) node[below] {\small{3}};
\draw[->, thin] (1,2.5) .. controls +(up:0.5cm)  .. (2,2.5);
\fill (2,2.5) circle (2pt) node[below] {\small{4}};
\fill (3,2.5) circle (2pt) node[below] {\small{2}};
\fill (4,2.5) circle (2pt) node[below] {\small{3}};
\draw[->, thin] (4,2.5) .. controls +(up:0.5cm)  .. (7,2.5);
\fill (5,2.5) circle (2pt) node[below] {\small{1}};
\draw[-,very thick] (5.3,2.3) -- (5.3,2.7);
\draw[->, thin] (5,2.5) .. controls +(up:0.5cm)  .. (8,2.5);
\fill (6,2.5) circle (2pt) node[below] {\small{3}};
\fill (7,2.5) circle (2pt) node[below] {\small{$\not 2\,\bar 3$}};
\fill (8,2.5) circle (2pt) node[below] {\small{$\not 2\,\bar 1$}};
\fill (9,2.5) circle (2pt) node[below] {\small{1}};
\fill (10,2.5) circle (2pt) node[below] {\small{4}};
\fill (11,2.5) node[right]{\small{$W_6 =31,\quad B_6=31,\quad F_6=3\bar 3\bar 1$}};
 %
\draw[-] (-1,1) -- (11,1);
\fill (0,1) circle (2pt) node[below] {\small{1}};
\draw[->, thin] (0,1) .. controls +(up:0.5cm)  .. (3,1);
\fill (1,1) circle (2pt) node[below] {\small{3}};
\draw[->, thin] (1,1) .. controls +(up:0.5cm)  .. (2,1);
\fill (2,1) circle (2pt) node[below] {\small{4}};
\fill (3,1) circle (2pt) node[below] {\small{2}};
\fill (4,1) circle (2pt) node[below] {\small{3}};
\draw[->, thin] (4,1) .. controls +(up:0.5cm)  .. (7,1);
\fill (5,1) circle (2pt) node[below] {\small{1}};
\draw[->, thin] (5,1) .. controls +(up:0.5cm)  .. (8,1);
\fill (6,1) circle (2pt) node[below] {\small{3}};
\draw[-,very thick] (6.3,0.8) -- (6.3,1.2);
\draw[->, thin] (6,1) .. controls +(up:0.5cm)  .. (10,1);
\fill (7,1) circle (2pt) node[below] {\small{$\not 2\,\bar 3$}};
\fill (8,1) circle (2pt) node[below] {\small{$\not 2\,\bar 1$}};
\fill (9,1) circle (2pt) node[below] {\small{1}};
\fill (10,1) circle (2pt) node[below] {\small{$\not 4\,\bar 3$}};
\fill (11,1) node[right] {\small{$W_7=313,\quad B_7=313,\quad F_6=\bar 3\bar 11\bar 3$}};
 %
\draw[-] (-1,-0.5) -- (11,-0.5);
\fill (0,-0.5) circle (2pt) node[below] {\small{1}};
\draw[->, thin] (0,-0.5) .. controls +(up:0.5cm)  .. (3,-0.5);
\fill (1,-0.5) circle (2pt) node[below] {\small{3}};
\draw[->, thin] (1,-0.5) .. controls +(up:0.5cm)  .. (2,-0.5);
\fill (2,-0.5) circle (2pt) node[below] {\small{4}};
\fill (3,-0.5) circle (2pt) node[below] {\small{2}};
\fill (4,-0.5) circle (2pt) node[below] {\small{3}};
\draw[->, thin] (4,-0.5) .. controls +(up:0.5cm)  .. (7,-0.5);
\fill (5,-0.5) circle (2pt) node[below] {\small{1}};
\draw[->, thin] (5,-0.5) .. controls +(up:0.5cm)  .. (8,-0.5);
\fill (6,-0.5) circle (2pt) node[below] {\small{3}};
\draw[->, thin] (6,-0.5) .. controls +(up:0.5cm)  .. (10,-0.5);
\fill (7,-0.5) circle (2pt) node[below] {\small{$\not 2\,\bar 3$}};
\draw[-,very thick] (7.3,-0.7) -- (7.3,-0.3);
\fill (8,-0.5) circle (2pt) node[below] {\small{$\not 2\,\bar 1$}};
\fill (9,-0.5) circle (2pt) node[below] {\small{1}};
\fill (10,-0.5) circle (2pt) node[below] {\small{$\not 4\,\bar 3$}};
\fill (11,-0.5) node[right]{\small{$W_8=13,\quad B_7=13\bar 3,\quad F_6=\bar 11\bar 3$}};
 %
\end{tikzpicture}
\caption[smallcaption]{An arrival scenario on the graph of Figure \ref{fig:example1}, 
and the trajectories of the 
three Markov chains.} 
\label{fig:example2}
\end{center}
\end{figure} 
\end{ex}

\subsection{Reversibility}
\label{subsec:reverse}
 For both chains $\suite{B_n}$ and $\suite{F_n}$, a state $\bw \in \bV^*$ 
is said admissible if it can be reached by the chain under consideration under {\sc fcfm}. We denote 
\begin{align}
\mbB &:=\Bigl\{\bw \in \bV^*\,:\,\bw \mbox{ is admissible for }\suite{B_n}\Bigl\}\label{eq:admB};\\
\mbF &:=\Bigl\{\bw \in \bV^*\,:\,\bw \mbox{ is admissible for }\suite{F_n}\Bigl\}.\nonumber
\end{align}
Define the following measure on $\bV^*$, 
\begin{equation}
\left\{\begin{array}{ll}
\nu_B(\emptyset) &=1;\\
\nu_B\left(\bw\right) &=\prod\limits_{i=1}^p \mu(i)^{|\bw|_i+| \td \bw|_i},\,\bw \in \bV^*\setminus\{\emptyset\}. 
\end{array}\right.
\label{eq:PiB}
\end{equation}
Observe that by the very definition (\ref{eq:PiB}), the measure of a word $\bw$ does not change whenever any of its letters $a$ is exchanged with $\td a$. 
In particular, we have $\nu_B(\cv{\td \bw})=\nu_B(\bw)$ for any $\bw$ in $\bV^*$. We have the following result, 
\begin{proposition}
\label{prop:product}
Suppose that $\mu$ belongs to {\sc Ncond}$(G)$ defined by (\ref{eq:Ncond}). 
Then the measure $\nu_B$ defined by (\ref{eq:PiB}) is finite on $\mbB$. Moreover, the Backwards detailed Markov chain $\suite{B_n}$ 
is positive recurrent, and admits $\beta\nu_B$ as a stationary measure for any $\beta>0$. 
\end{proposition}
 
Before showing Proposition \ref{prop:product} we need to introduce a few technical results.
\begin{lemma}
\label{lemma:M83}
Let $\bw=\bw_1...\bw_q\in \bV^*$. Then $\bw\in \mbB$ if, and only if the following three conditions hold, 
\begin{align}
\label{eq:admissible-1}
&\bw_{1}\in \maV;\\ 
\label{eq:admissible0}
&\bw_{k} \pv \bw_{\ell}\mbox{ for any }k,\ell \in \llbracket 1,q \rrbracket\mbox{ such that } \bw_{k}\in \maV\mbox{ and }\bw_\ell\in \maV;\\
\label{eq:admissible}
&\bw_{k} \pv \td{\bw_{j}}\mbox{ for any } k,j,\,\, 1\le k <j \le q,\mbox{ such that } \bw_{k}\in \maV\mbox{ and }\td{\bw_{j}}\in \td\maV. 
\end{align}  
\end{lemma}
\begin{proof}[Proof of Lemma \ref{lemma:M83}] 
By the very definition of $\suite{B_n}$, (\ref{eq:admissible-1}) holds true for any (non-empty) admissible word of the chain. Moreover the necessity of (\ref{eq:admissible0}) is obvious: would an element $\bw$ of $\mbB$ contain two compatible letters in $\maV$, the two corresponding items would have been matched. 
Let us prove the necessity of (\ref{eq:admissible}). 
Fix two such indexes $k$ and $j$ in $\llbracket 1,q \rrbracket$. 
This means that $V_{i(n)+j-1}$ is matched with an item $V_\ell$ of class $\td{\bw_{j}}$. 
Then, it can be easily checked that 
$\bw_{k} \v \td{\bw_{j}}$ implies that $V_{\ell}$ is matched with $V_{i(n)+k-1}$, an absurdity since $V_{i(n)+k-1}$ of class $\bw_{k}$ is still unmatched at $n$. 
This completes the proof of necessity. 

Regarding sufficiency, fix a non-empty state $\bw$ satisfying (\ref{eq:admissible-1}), (\ref{eq:admissible0}) and (\ref{eq:admissible}):  
\[\bw = b_1\,\td{a_{11}}\,\td{a_{12}}\,...\,\td{a_{1k_1}}\,b_2\,\td{a_{21}}...\td{a_{2k_2}}\,b_3....b_q\,\td{a_{q1}}...\td{a_{qk_q}},\]
where $q\ge 1$, $k_q \in \N$ for all $q$, $b_\ell \in \maV$ for all $\ell$ and $a_{\ell j} \in \maV$ for all $\ell,j$. 
In particular, from (\ref{eq:admissible0}) we have that $b_i \pv b_j$ for any $i \ne j$ whereas from (\ref{eq:admissible}), $a_{\ell j} \pv b_i$ for any $j$ and any $i \le \ell$. 
Let us show that the chain $\suite{B_n}$ can reach the state $\bw$. For this we construct inductively an arrival vector $A$ leading to $\bw$ from the state $\emptyset$. 
At first, set 
\[A:=\left(V_{n-q+1},V_{n-q+2},...,V_{n}\right)=\left(b_1,...,b_q\right).\]
Then, we investigate all elements $\td{a_{\ell j}}$ from left to right, as follows. We start from $a_{11}$: 

(1) If there is no element $\td{a_{\ell' j'}}$ to the right of $\td{a_{11}}$ such that $a_{11} \v a_{\ell' j'}$, then set 
$V_{n-q-1}=a_{11}$, $V_{n-q}=b_1$, $V_{n-q+1}=c \in \maE\left(a_{11}\right)$, in a way that the $a_{11}$-item is matched with that of class $b$, and 
we retrieve the letter $\td{a_{11}}$ to the right of $b_1$ because the $c$-item is matched with $a_{11}$;  

(2) If there exists an element $\td{a_{\ell' j'}}$ to the right of $\td{a_{11}}$ such that $a_{11} \v a_{\ell' j'}$, 
we investigate all terms $\td{a_{\ell'' j''}}$ to the left of $\td{a_{\ell' j'}}$. If one of them, $\td{a_{\ell'' j''}}$ is such that 
$a_{\ell'' j''} \v a_{\ell' j'}$, then $a_{\ell' j'}$ can be matched in FCFM with another item than $a_{11}$. Then we do as in case (1): 
$V_{n-q-1}=a_{11}$, $V_{n-q}=b_1$, $V_{n-q+1}=c \in \maE\left(a_{11}\right)$; 

(3) If an element $\td{a_{\ell' j'}}$ to the right of $\td{a_{11}}$ is such that $a_{11} \v a_{\ell' j'}$, 
and no term $\td{a_{\ell'' j''}}$ to the left of $\td{a_{\ell',j'}}$ is such that 
$a_{\ell'' j''} \v a_{\ell'j'}$, we interpose a term $a_{\ell' j'}$ in $A$ between $b_1$ and $b_2$ 
(or at the extreme right of $A$ if $q=1$), and a term $a_{11}$ between 
$b_{\ell'}$ and $b_{\ell'+1}$ (or at the extreme right of $A$ if $q=\ell'$), in a way that the two corresponding items are matched. 

Then, by induction we can investigate in the same way all the terms $\td{a_{\ell j}}$ not yet considered, 

(1) If there is no element $\td{a_{\ell'  j'}}$ to the right of $\td{a_{\ell j}}$ such that $a_{\ell j} \v a_{\ell' j'}$, then in $A$ 
we interpose a letter $a_{\ell j}$ just to the left of $b_1$, and a letter $c \in \maE\left(a_{\ell j}\right)$ to the left of $b_{\ell+1}$ if $\ell <q$ 
(or at the extreme right of $A$ if $\ell=q$); the items of classes $c$ and $a_{\ell,j}$ are matched and a term $\td{a_{\ell j}}$ appears at the right place 
in the detailed state of the system; 

(2) We do as in case (1) if there exists an element $\td{a_{\ell' j'}}$ to the right of $\td{a_{\ell j}}$ such that $a_{\ell j} \v a_{\ell' j'}$, 
but one of the non yet investigated terms $\td{a_{\ell' j'}}$ to the left of $\td{a_{\ell'j'}}$ is such that 
$a_{\ell'' j''} \v a_{\ell'j'}$. 

(3) If there exists an element $\td{a_{\ell' j'}}$ to the right of $\td{a_{\ell j}}$ such that $a_{\ell j} \v a_{\ell' j'}$, 
and no term $\td{a_{\ell'' j''}}$ to the left of $\td{a_{\ell'j'}}$ is such that 
$a_{\ell'' j''} \v a_{\ell'j'}$, then in $A$ we interpose a term $a_{\ell' j'}$ just to the left of $b_{\ell +1}$ (or at the extreme right of $A$ if 
$\ell =q$), and a term $a_{\ell j}$ to the immediate left of $b_{\ell'+1}$ (or at the extreme right of $A$ if 
$\ell' =q$), in a way that the two corresponding items are matched. 

We continue this construction, 
until all the letters $\td{a_{\ell j}} \in \bw|_{\td \maV}$ are investigated and the corresponding items are matched. 
The resulting vector $A$ that we obtain is of size $q'\in \left\llbracket q+\sum_{\ell=1}^q k_{\ell},q+2\sum_{\ell=1}^q k_{\ell}\right\rrbracket$. 
Indeed, the number of 
items added to $A$ is at least equal to the number of letters of $\bw|_{\td \maV}$, and at most equal to twice the latter number (which is the case 
if all the corresponding items entered the system before the item of class $b_1$ and are matched after the arrival time of the latter). 
Finally, for any $n \ge q+2\sum_{\ell=1}^q k_{\ell}$, if $B_{n-q'}=\emptyset$, an arrival scenario $A$ for the $q'$ following time epochs yields to a state $B_n=\bw$. This concludes the proof.  
\end{proof} 

As a consequence, 
\begin{lemma}
\label{lemma:admissible}
The mapping $\bw \mapsto \cv{\td \bw}$ is one-to-one from $\mbB$ into $\mbF$.
%
%
\end{lemma}

\begin{proof}
From Lemma \ref{lemma:M83}, it is sufficient to prove that a state $\bw$ belongs to $\mbF$ if and only if 
$\cv{\td{\bw}}$ satisfies (\ref{eq:admissible-1}), (\ref{eq:admissible0}) and (\ref{eq:admissible}). The proof of this statement is similar to that of Lemma \ref{lemma:M83}.  
\end{proof} 

We can now state the following connection between the dynamics of $\suite{B_n}$ and $\suite{F_n}$, 
\begin{lemma}
\label{lemma:kelly}
Let $\nu_B$ be the measure on $\maV^*$ defined by (\ref{eq:PiB}). Then,
for any $\bw,\bw' \in\mbB$ such that $\pr{B_{n+1}=\bw' | B_n=\bw} >0$ 
we have that 
\begin{equation}
\label{eq:kelly}
\nu_B(\bw)\pr{B_{n+1}=\bw' | B_n=\bw} = \nu_B\left(\cv{\td{\bw'}}\right)\pr{F_{n+1}=\cv{{\td \bw}} | F_n=\cv{\td{\bw'}}}.
\end{equation}
\end{lemma}

\begin{proof} 
Fix $\bw \in \mbB$ (so that $\cv{\td{\bw}} \in \mbF$ from Lemma \ref{lemma:admissible}). 
Whenever $\bw\ne \emptyset$ we set 
$\bw=\bw_1...\bw_q\in \bV^*$, and recall that $\bw_1 \in \maV$ in that case. 
There are five possible cases for the transition of $\suite{B_n}$, which we address one by one, 

\medskip

(1) Suppose that $\bw\ne \emptyset$ and $\bw'=\bw a$ for some $a\in \maV$. Plainly, such a state is admissible if and only if 
$a \in \maE\left(\maV(\bw)\right)^c$. The backwards chain moves from $\bw$ to $\bw a$ at $n+1$ whenever $V_{n+1}=a$, so we have that $\pr{B_{n+1}=\bw a | B_n=\bw}=\mu(a).$ 
On the other hand, $F_n=\cv{{\td{{\bw} a}}}$ entails that the item entering at $n+1$ is matched with an item of class $a$ entered up to time $n$. 
Therefore we necessarily have that $F_{n+1}=\cv{\td\bw}$, in other words 
\begin{align*}
\nu_B\left(\cv{\td{{\bw}a}}\right)\pr{F_{n+1}=\cv{\td \bw} | F_n=\cv{\td{{\bw}a}}}&=\nu_B\left(\cv{\td{{\bw}a}}\right)\\
                                                                   &=\nu_B(\bw)\mu(a)
																   =\nu_B(\bw)\pr{B_{n+1}={\bw}a | B_n=\bw}.
																   \end{align*}

(2) Suppose now that $\bw\ne \emptyset$ and $\bw'=\bw_1...\bw_{k-1}\,\td a \,\bw_{k+1}...\bw_q\td{\bw_{k}}$.
This means that $\bw_{k}\in \maV$ and that the item $V_{n+1}$ is of class $a$, where 
$a \in\maE\left(\bw_{k}\right)\,\cap\, \maE\left(\maV\left(\bw_1....\bw_{k-1}\right)\right)^c,$ 
so that in FCFM, $V_{n+1}$ is matched with the item $V_{i(n)+k-1}$ of class $\bw_{k}$. 
Suppose that \[F_n=\cv{\td {\bw'}}=\bw_{k}\td{\bw_q}\,\td{\bw_{q-1}}\,...\td{\bw_{k+1}}\, a \,\td{\bw_{k-1}}\,...\,\td{\bw_1}.\] 
Then from Lemma \ref{lemma:M83}, $\bw_{k}$ is not adjacent to any of the elements of $\maV\left(\td{\bw_q}\,\td{\bw_{q-1}}\,...\,\td{\bw_{k+1}}\right)$. 
But $\bw_{k} \v a$, so the item $V_{n+1}$ of class $\bw_{k}$ is matched with the item $V_{n+q-k+2}$ of class $a$, and we have with probability 1,  
\[F_{n+1}= \td{\bw_q}\,\td{\bw_{q-1}}\,...\td{\bw_{k+1}}\, \td{\bw_{k}} \,\td{\bw_{k-1}}\,...\,\td{\bw_1}=\cv{\td \bw}.\]
 Therefore, in this case, 
\begin{align*}
\nu_B\left(\cv{\td {\bw'}}\right)\pr{F_{n+1}=\cv{\td \bw} | F_n=\cv{\td {\bw'}}} &=\nu_B\left(\cv{\td{\bw'}}\right)\\
                                                                   &=\nu_B\left(\cv{\td \bw}\right)\mu(a)\\
                                                                   &=\nu_B(\bw)\mu(a)
																   =\nu_B(\bw)\pr{B_{n+1}=\bw' | B_n=\bw}.
																   \end{align*} 
(3) Now, suppose that $\bw\ne \emptyset$ and $\bw'=\bw_{k}\bw_{k+1}...\bw_q\td{\bw_1}$ for some $k \in \llbracket 2,q \rrbracket$. 
This means that the class of item $V_{n+1}$ belongs to $\maE\left(\bw_1\right)$, so $V_{n+1}$ is matched with the oldest item 
in line $V_{i(n)}$. Then $\bw_2,\bw_{3},...,\bw_{k-1}$ all belong to $\td{\maV}$, and so $\bw_{k} \in \maV$ and is the class of the oldest item in line after $V_{i(n)}$, now becoming the new oldest one. 
Suppose that \[F_n=\cv{\td {\bw'}}=\bw_1\td{\bw_q}\,\td{\bw_{q-1}}\,...\td{\bw_{k+1}}\,\td{\bw_{k}}.\] Applying again Lemma \ref{lemma:M83}, we obtain that 
$\bw_1$ is not adjacent to any of the elements of the set 
$\maV\left(\td{\bw_q}\,\td{\bw_{q-1}}\,...\,\td{\bw_{k}}\right)$, so the state $\cv{\td{\bw'}}$ clearly belongs to $\mbF$. All the same, again in view of Lemma \ref{lemma:M83}, $\bw_1$ is not adjacent to any of the elements $\td{\bw_2},\td{\bw_{3}},....,\td{\bw_{k-1}}$, which all are elements of $\maV$. So we obtain 
\[F_{n+1}=\td{\bw_q}\,\td{\bw_{q-1}}\,...\td{\bw_{k+1}}\,\td{\bw_k}\,\td{\bw_{k-1}}\,...\td{\bw_2}\,\td{\bw_1}=\cv{\td\bw}\]
if the incoming items $V_{n+q-k+3}$,...,$V_{n+q}$ are of respective classes $\td{\bw}_{k-1},...,\td{\bw}_{2}$ and $V_{n+q+1}$ is of a class 
belonging to $\maE\left(\bw_1\right)$. This occurs with probability 
\[\mu\Bigl(\td{\bw_{k-1}}\Bigl)...\mu\Bigl(\td{\bw_2}\Bigl)\mu\Bigl(\maE\left(\bw_1\right)\Bigl).\]
Gathering all the above we obtain that
\begin{align*}
\nu_B\left(\cv{\td {\bw'}}\right)\pr{F_{n+1}=\cv{\td \bw} | F_n=\cv{\td {\bw'}}} &=\nu_B\left(\cv{\td {\bw'}}\right)\mu\Bigl(\td{\bw_{k-1}}\Bigl)...\mu\Bigl(\td{\bw_2}\Bigl)
\mu\Bigl(\maE\left(\bw_1\right)\Bigl)\\
                                                                   &=\nu_B\left(\cv{\td \bw}\right)\mu\Bigl(\maE\left(\bw_1\right)\Bigl)\\
                                                                   &=\nu_B(\bw)\mu\Bigl(\maE\left(\bw_1\right)\Bigl)
																   =\nu_B(\bw)\pr{B_{n+1}=\bw' | B_n=\bw}.
																   \end{align*}																 
(4) Suppose now that $\bw\ne \emptyset$ and $\bw'=\emptyset$, which is possible only if $\bw_1\in \maV$, 
the incoming item $V_{n+1}$ belongs to $\maE\left(\bw_1\right)$, and if $q \ge 2$,  
$\bw_2,...,\bw_q \in \td\maV$, which implies again from Lemma \ref{lemma:M83} that $\bw_1 \pv \td{\bw_j}$ for any $j\in \llbracket 2,q \rrbracket$. 
Thus, $F_n=\emptyset$ leads to the state $F_{n+1}=\cv{\td {\bw}}=\td{\bw_q}\,\td{\bw_{q-1}}\,....\,\td{\bw_2}\,\td{\bw_1}$ 
if and only if $V_{n+1}$ is of class $\bw_1$, and then $V_{n+2}$ is of class $\td{\bw_q}$, $V_{n+3}$ is of class $\td{\bw_{q-1}}$, 
and so on ..., $V_{n+q}$ is of class $\td{\bw_2}$ and $V_{n+q+1}$ is of a class belonging to $\maE\left(\bw_1\right).$ 
This event occurs with probability $\mu(\bw_1)\mu\left(\td{\bw_q}\right)....\mu\left(\td{\bw_2}\right)\mu\left(\maE\left(\bw_1\right)\right).$  
So we obtain   
\begin{align*}
\nu_B\left(\emptyset\right)\pr{F_{n+1}=\cv{\td \bw} | F_n=\emptyset} &=\mu\left(\bw_1\right)\mu\Bigl(\td{\bw_2}\Bigl)\mu\Bigl(\td{\bw_3}\Bigl)...\mu\Bigl(\td{\bw_q}\Bigl)\mu\Bigl(\maE\left(\bw_1\right)\Bigl)\\
                                                                   &=\nu_B(\bw)\mu\Bigl(\maE\left(\bw_1\right)\Bigl)
																   =\nu_B(\bw)\pr{B_{n+1}=\bw' | B_n=\bw}.
																   \end{align*}																   

(5) The only case that remains to be treated is when $\bw=\emptyset$. 
Then for any $a\in \maV$, we obtain $B_{n+1}=a$ provided that $V_{n+1}$ is of class $a$, which occurs 
with probability $\mu(a)$.  
Then, $F_n=\td{a}$ means that $V_{n+1}$ is matched with an item of class $a$ that was entered before $n$. 
Then we necessarily have that $F_{n+1}=\emptyset$, and thus  
\[
\nu_B\left(\td{a}\right)\pr{F_{n+1}=\emptyset | F_n=\td{a}}=\nu_B\left(a\right)=\mu(a) =\nu_B(\emptyset)\pr{B_{n+1}=a | B_n=\emptyset}.   \]
This completes the proof. 
%
%
%
\end{proof}

We can now turn to the proof of Proposition \ref{prop:product}. 
\begin{proof}[Proof of Proposition \ref{prop:product}]
Suppose that $\mu$ belongs to the set $\textsc{Ncond}(G)$ defined by 
(\ref{eq:Ncond}). From Lemma \ref{lemma:M83}, we know that for any word $\bw$ in $\mbB$, $\maV(\bw)$ is an element of $\mathbb I(G)$, that is, the letters of $\maV$ present in $\bw$ form an independent set of $\maV$. Also, the intermediate letters of $\bw$ in $\bar{\maV}$ have counterparts in $\maV$ that are not adjacent to any prior letter of 
$\bw$ in $\maV$. 
Therefore, we can partition the set $\mbB$ as follows: Denote for any independent set 
$\maI:=\left\{i_1,...,i_{|\maI|}\right\}\in \mathbb I(G)$, and for any permutation $\sigma$ in $\mathfrak S_{|\maI|}$, 
by $\mbB_{\maI,\sigma}$, the set of non-empty words $\mathbf w\in \mbB$ whose letters in $\maV$ are exactly the elements of $\maI$ (i.e., $\maV(\mathbf w)=\maI$), and 
$i_{\sigma(1)}$, $i_{\sigma(2)}$, ... $i_{\sigma(|\maI|)}$ denote the letters of $\maV(\mathbf w)$ ranked in their order of first occurence in the word $\mathbf w$. In other words, 
in $\bw$ read from left to right, the letter $i_{\sigma(1)}$ appears first (possibly several times), the first letter in $\maV(\mathbf w)\setminus\left\{i_{\sigma(1)}\right\}$ is 
$i_{\sigma(2)}$, the first letter in $\maV(\mathbf w)\setminus\left\{i_{\sigma(1)},i_{\sigma(2)}\right\}$ is $i_{\sigma(3)}$, and so on. 
Then, we clearly have the following disjoint union, 
\begin{equation}
\label{eq:partitionB}
\mathbb B = \{\emptyset\} \cup \bigcup_{\maI\in\mathbb I(G)} \bigcup_{\sigma\in\mathfrak S_{|\maI|}} \mathbb B_{\maI,\sigma}.
\end{equation}
Now, fix $\maI\in\mathbb I(G)$ and $\sigma\in\mathfrak S_{|\maI|}$. Then, for any word $\bw$ in $\mbB_{\maI,\sigma}$, and any $j\in\llbracket 1,|\maI|\rrbracket$, between the first occurrence of the letter $i_{\sigma(j)}$ (excluded) and the first occurrence of the letter $i_{\sigma(j+1)}$ (or to the right of $i_{\sigma(|\maI|)}$ if $j\equiv |\maI|$), there are, say, $\varphi(j)$ letters in $\bw$ (where $\varphi(j) \in\N$) that belong to the set 
$\{i_{\sigma(1)},...,i_{\sigma(j)}\}$. Denote these letters (again, read from left to right) by $k_{1},...,k_{\varphi(j)}$, and set $k_0:=i_{\sigma(j)}$. 
In between these letters, say between $k_l$ and $k_{l+1}$ for $\varphi(j)>0$ and $l\in\llbracket 0,\varphi(j)-1 \rrbracket$ 
(or to the right of $k_{\varphi(j)}$ for $l \equiv \varphi(j)$), 
there can be in $\bw$, any number $\psi(l) \in\N$ of letters in $\bar{\maV}$ whose counterpart in $\maV$ are element of 
$\maE\left(\left\{i_{\sigma(1)},...,i_{\sigma(j)}\right\}\right)^c$. 
Denote by $\left\{h_1,...,h_{r(j)}\right\}$, the elements of $\maE\left(\left\{i_{\sigma(1)},...,i_{\sigma(j)}\right\}\right)^c$, whenever the latter is non empty. 
Denoting also for all $m=1,...,r(j)$, by $p_m,$ the number of occurences of the letter $\td{h_m}$ in between $k_l$ and $k_{l+1}$ in the word $\bw$, we obtain the following, 
\begin{align*}
\nu_B\left(\mbB_{\maI,\sigma}\right)
&=\prod_{j=1}^{|\maI|} \sum_{\varphi(j)\in\N}\sum_{\substack{(k_1,...,k_{\varphi(j)})\\\in \{i_{\sigma(1)},...,i_{\sigma(j)}\}^{\varphi(j)}}}\prod_{l=0}^{\varphi(j)} \mu(k_l) \sum_{\psi(l)\in\N} 
\displaystyle\sum_{\substack{p_1,\,...,\,p_{r(j)}\in\N:\\p_1+...+p_{r(j)} = \psi(l)}} {\psi(l)\choose p_1,...,p_{r(j)}} \displaystyle\prod_{m=1}^{r(j)}\mu\left(h_m\right)^{p_m} \nonumber\\
&=\prod_{j=1}^{|\maI|} \sum_{\varphi(j)\in\N}\sum_{\substack{(k_1,...,k_{\varphi(j)})\\\in \{i_{\sigma(1)},...,i_{\sigma(j)}\}^{\varphi(j)}}}\prod_{l=0}^{\varphi(j)} \mu(k_l) \sum_{\psi(l)\in\N} 
\left(\sum_{m=1}^{r(j)} \mu(h_m)\right)^{\psi(l)} \nonumber\\
&= \prod_{j=1}^{|\maI|} \sum_{\varphi(j)\in\N}\sum_{\substack{(k_1,...,k_{\varphi(j)})\\\in \{i_{\sigma(1)},...,i_{\sigma(j)}\}^{\varphi(j)}}}\prod_{l=0}^{\varphi(j)} {\mu(k_l) \over \mu\left(\maE\left(\{i_{\sigma(1)},...,i_{\sigma(j)}\}\right)\right)},
\end{align*}
where products and sums over empty sets are respectively understood as 1 and 0. 
Now, if $\varphi(j)>0$, denoting for all $n=1,...,j$, by $q_n,$ the number of occurences of the letter 
$i_{\sigma(n)}$ in the word $k_1...k_{\varphi(j)}$, the above becomes 
\begin{align*}
\nu_B\left(\mbB_{\maI,\sigma}\right)
&=\prod_{j=1}^{|\maI|} 
{\mu\left(i_{\sigma(j)}\right) \over \mu\left(\maE\left(\left\{i_{\sigma(1)},...,i_{\sigma(j)}\right\}\right)\right)}\displaystyle\sum_{\varphi(j)\in\N}
{ \displaystyle\sum_{\substack{q_1,\,...,\,q_j\in\N:\\q_1+...+q_j = \varphi(j)}} {\varphi(j)\choose q_1,...,q_j} \displaystyle\prod_{n=1}^j\mu\left(i_{\sigma(n)}\right)^{q_n} \over \mu\left(\maE\left(\{i_{\sigma(1)},...,i_{\sigma(j)}\}\right)\right)^{\varphi(j)}}
\nonumber\\
&=  \prod_{j=1}^{|\maI|} 
{\mu\left(i_{\sigma(j)}\right) \over \mu\left(\maE\left(\left\{i_{\sigma(1)},...,i_{\sigma(j)}\right\}\right)\right)}
\sum_{\varphi(j)\in\N} \left({\mu\left(\left\{i_{\sigma(1)},...,i_{\sigma(j)}\right\}\right) \over \mu\left(\maE\left(\{i_{\sigma(1)},...,i_{\sigma(j)}\}\right)\right)}\right)^{\varphi(j)}
\nonumber\\
&= \prod_{j=1}^{|\maI|} 
{\mu\left(i_{\sigma(j)}\right)\over \mu\left(\maE\left(\left\{i_{\sigma(1)},...,i_{\sigma(j)}\right\}\right)\right)}
{\mu\left(\maE\left(\{i_{\sigma(1)},...,i_{\sigma(j)}\}\right)\right) \over \mu\left(\maE\left(\{i_{\sigma(1)},...,i_{\sigma(j)}\}\right)\right)-\mu\left(\left\{i_{\sigma(1)},...,i_{\sigma(j)}\right\}\right)}
\nonumber\\
&= \prod_{j=1}^{|\maI|} 
{\mu\left(i_{\sigma(j)}\right) \over \mu\left(\maE\left(\{i_{\sigma(1)},...,i_{\sigma(j)}\}\right)\right)-\mu\left(\left\{i_{\sigma(1)},...,i_{\sigma(j)}\right\}\right)},
\end{align*}
where the assumption that  $\mu\in\textsc{Ncond}(G)$ is used in the third equality, as $\left\{i_{\sigma(1)},...,i_{\sigma(j)}\right\}$ is an independent set of $G$ for any $j$ as above. 
As a conclusion, from (\ref{eq:partitionB}) we obtain that 
\begin{equation}
\label{eq:alpha0}
\nu_B(\mbB)= 1 +\sum_{\maI\in\mathbb I(G)} \sum_{\sigma\in\mathfrak S_{|\maI|}} \prod_{j=1}^{|\maI|} 
{\mu\left(i_{\sigma(j)}\right) \over \mu\left(\maE\left(\{i_{\sigma(1)},...,i_{\sigma(j)}\}\right)\right)-\mu\left(\left\{i_{\sigma(1)},...,i_{\sigma(j)}\right\}\right)},
\end{equation}
and so $\nu_B$ is a finite measure on $\mbB$. 

Then, it suffices to apply Kelly's Lemma (\cite{kelly:79}, Section 1.7): define for any 
$\bw,\bw'\in \mbB$,
\begin{equation}
\label{eq:defP}
P_{\bw',\bw}={\pr{B_{n+1}=\bw' \mid B_n=\bw}\nu_B(\bw) \over \nu_B(\bw')}.
\end{equation}
Then, $\nu_B$ is a stationary distribution of $\suite{B_n}$ if $P$ defines a transition operator on $\mbB$. But this is a simple consequence of 
Lemmas \ref{lemma:admissible} and \ref{lemma:kelly}: for any $\bw' \in \mbB$ we have that 
\begin{align*}
\sum_{\bw \in \mbB} P_{\bw',\bw} 
                            &= \sum_{\bw \in\mbB} {\nu_B\left(\cv{\td{\bw'}}\right)\pr{F_{n+1}=\cv{{\td \bw}} | F_n=\cv{\td{\bw'}}} \over \nu_B(\bw')}\\
                            &= \sum_{\bw \in \mbB} \pr{F_{n+1}=\cv{{\td \bw}} | F_n=\cv{\td{\bw'}}}
                            =1,
\end{align*}
which concludes the proof
\end{proof}

\begin{remark}\rm 
As is easily seen, we can exchange the roles of $B_n,B_{n+1}$ and $F_n,F_{n+1}$ in (\ref{eq:defP}), implying that $\suite{F_n}$ is the reversed Markov chain of $\suite{B_n}$, on a sample space where arrivals are reversed 
in time and exchanged with their match. In particular, the forward chain $\suite{F_n}$ also admits $\nu_B$ as a stationary measure on $\mbF$. 
\end{remark}

We are now in position to prove our main result.  
\subsection{Proof of Theorem \ref{thm:main}}
\label{subsec:stabFCFM}
We know from Proposition 2 of \cite{MaiMoy16} that
if $\mu$ is not an element of {\sc Ncond}($G$), then the chain $\suite{W_n}$ is transient or null recurrent. 
If we now assume that $\mu\in\textsc{Ncond}(G)$, then, first, observe that 
in view of (\ref{eq:alpha0}), the measure defined for all $\bw\in\mbB$ by $\alpha\nu_B(\bw)$, for $\alpha$ in (\ref{eq:alpha}) and $\nu_B$ defined by (\ref{eq:PiB}), defines a probability measure on $\mbB$. 
Second, from Proposition \ref{prop:product}, the auxiliary chain 
$\suite{B_n}$ is positive recurrent, and admits $\nu_B$ as a stationary measure. So from (\ref{eq:WB}), $\suite{W_n}$ is also positive recurrent. 
As it is also irreducible on $\mathbb W$, it has a unique stationary probability. To check that the latter is given by $\Pi_W$ defined by (\ref{eq:PiW}), it thus suffices to check that 
\begin{equation*}
\Pi_{W}(w)=\alpha\sum\limits_{\bw\in \mbB:\,\bw|_{\maV}=w} \nu_B(\bw)\,\quad\mbox{for any }w\in \mathbb W.
\end{equation*} 
Let $w=w_1...w_q \in \mathbb W$. From Lemma \ref{lemma:M83}, any $\bw \in \mbB$ such that $\bw|_{\maV}=w$ is of the form 
\[\bw = w_1\td{a_{\small{11}}}\,\td{a_{\small{12}}}\,...\,\td{a_{\small{1k_1}}}\,w_2\,\td{a_{\small{21}}}...\td{a_{\small{2k_2}}}\,w_3\,...\,w_q\td{a_{\small{q1}}}\,...\,\td{a_{\small{qk_q}}},\]
where any of the elements $a_{\small{\ell j}}$ is such that $a_{\ell j} \pv w_{i}$ for any $i \le \ell$.  
We then obtain that 
\begin{align*}
\alpha\sum\limits_{\bw\in \mbB:\,\bw|_{\maV}=w} \nu_B(\bw) 
&= \alpha\prod\limits_{\ell =1}^q \left(\mu(w_\ell)\left(1+\sum_{k\in\N_+}\sum\limits_{(a_{\ell 1},...,a_{\ell k}) \in \left(\maE\left(\left\{w_1,...,w_\ell\right\}\right)^c\right)^k}\prod_{j=1}^{k} \mu\left(a_{\ell j}\right)\right)\right)\\
&=\alpha\prod\limits_{\ell =1}^q \left(\mu(w_\ell)\sum_{k\in\N}\left(\mu\biggl(\maE\Bigl(\left\{w_1,...,w_\ell\right\}\Bigl)^c\biggl)\right)^{k}\right)\\
&=\alpha\prod\limits_{\ell =1}^q {\mu(w_\ell) \over 1-\mu\biggl(\maE\Bigl(\left\{w_1,...,w_\ell\right\}\Bigl)^c\biggl)}\\
&=\alpha\prod\limits_{\ell =1}^q {\mu(w_\ell) \over \mu\biggl(\maE\Bigl(\left\{w_1,...,w_\ell\right\}\Bigl)\biggl)}=\Pi_{W}(w). 
\end{align*}
\qed 

\subsection{Characteristics at equilibrium}
\label{subsec:perf}
We can easily deduce from Theorem \ref{thm:main}, closed form formulas for performance measures of the system in steady state. 
Denote by $W_\infty$, the stationary queue detail of the system, that is, a random variable distributed following the stationary probability $\Pi_W$. 

First, for any class $i\in\maV$, the average number of items of class $i$ in storage at equilibrium is given by 
\[\esp{\left|W_\infty\right|_i}= \sum_{k \in\N} k \Pi_W\left(\{w \in \mathbb W\,:\,|w|_i=k\}\right),\]
for $\Pi_W$ in (\ref{eq:PiW}). The average total number of items in storage at equilibrium then reads 
\[\esp{|W_\infty|}=\sum_{i\in\maV} \esp{\left|W_\infty\right|_i}= \sum_{k \in\N} k \Pi_W\left(\{w \in \mathbb W\,:\,|w|=k\}\right).\]

According to Little's law, for any $i\in\maV$, the waiting time before getting matched, for an item of class $i$ entering the system in steady state, is given by 
\[{\esp{|W_\infty|}\over \mu(i)} = {1\over \mu(i)} \sum_{k \in\N} k \Pi_W\left(\{w \in \mathbb W\,:\,|w|_i=k\}\right).\]
In particular, the probability that a class $i$-item does have to wait before getting matched in a stationary system reads 
\[\pr{|W_\infty|_{\maE(i)}=0} = \Pi_W\left(\{w \in \mathbb W\,:\,|w|_{\maE(i)}=0\}\right),\]
for $\Pi_W$ in (\ref{eq:PiW}).

\section{Perfect {\sc fcfm}-matchings}
\label{FCFMreverse}

\subsection{Perfect Infinite {\sc fcfm}-matchings}
\label{subsec:matching} 


As is easily seen, the model $(G,\mu,\textsc{fcfm})$ generates a family of random graphs, as follows. 
For any $n\in \N$, 
we consider the (initially empty) {\sc fcfm}-matching  
$$\mathbf M^{0}_{n}:=M^{\textsc{fcfm}}\left(V_{0}\,...\,V_{n-1}\right),$$ 
that is, the random graph in which the nodes are the items from $0$ to $n-1$, and there is an edge between two nodes if and only 
if the corresponding items are matched in {\sc fcfm}. The nodes of this random graph are labelled according to their classes, and (implicitly) to their arrival dates. 
However, for notational simplicity, we only keep the first labelling (i.e., the classes), and remove the second one, which is 
implicitly given by the order of nodes (from left to right) in Figures \ref{fig:increase} and \ref{Fig:exchange} below. 

The definition of a {\sc fcfm}-matching can then be extended to any time points $m<n$ in $\N$ such that $W_m=\emptyset$, by setting 
$$\mathbf M^{m}_{n}:=M^{\textsc{fcfm}}\left(V_{m}\,...\,V_{n-1}\right).$$
In $\mathbf M^{m}_{n}$, all nodes are 
of degree 0 or 1. The finite matching $\mathbf M^{m}_{n}$ is then said to be {\em perfect} 
if all of its nodes are of degree 1. 
It is then immediate that for any $m$ such that $W_m=\emptyset$ and any $n>m$, 
$\mathbf M^{m}_{n}$ is an induced subgraph of 
$\mathbf M^{m}_{n+1}$. See in Figure \ref{fig:increase}, an example concerning once again the graph $G$ of Figure \ref{fig:example1}. 

\begin{figure}[h!]
\begin{center}
\begin{tikzpicture}[scale=0.8]
%
\fill (-1,7) node{$\emptyset$};
\fill (5,7) node[right] {\small{$\mathbf M^0_0$}};
\fill (-1,6) circle (2pt) node[below] {\scriptsize{1}};
\fill (5,6) node[right] {\scriptsize{$\mathbf M^0_1$}};
\fill (-1,5) circle (2pt) node[below] {\scriptsize{1}};
\fill (0,5) circle (2pt) node[below] {\scriptsize{3}};
\fill (5,5) node[right] {\small{$\mathbf M^0_2$}};
\fill (-1,4) circle (2pt) node[below] {\scriptsize{1}};
\fill (0,4) circle (2pt) node[below] {\scriptsize{3}};
\fill (1,4) circle (2pt) node[below] {\scriptsize{4}};
\draw[->, thin] (1,4) .. controls +(up:0.5cm)  .. (0,4);
\fill (5,4) node[right] {\small{$\mathbf M^0_3$}};
\fill (-1,3) circle (2pt) node[below] {\scriptsize{1}};
\fill (0,3) circle (2pt) node[below] {\scriptsize{3}};
\fill (1,3) circle (2pt) node[below] {\scriptsize{4}};
\draw[->, thin] (1,3) .. controls +(up:0.5cm)  .. (0,3);
\fill (2,3) circle (2pt) node[below] {\scriptsize{2}};
\draw[->, thin] (2,3) .. controls +(up:0.5cm)  .. (-1,3);
\fill (5,3) node[right] {\small{$\mathbf M^0_4$}};
\fill (-1,2) circle (2pt) node[below] {\scriptsize{1}};
\fill (0,2) circle (2pt) node[below] {\scriptsize{3}};
\fill (1,2) circle (2pt) node[below] {\scriptsize{4}};
\draw[->, thin] (1,2) .. controls +(up:0.5cm)  .. (0,2);
\fill (2,2) circle (2pt) node[below] {\scriptsize{2}};
\draw[->, thin] (2,2) .. controls +(up:0.5cm)  .. (-1,2);
\fill (3,2) circle (2pt) node[below] {\scriptsize{3}};
\fill (5,2) node[right] {\small{$\mathbf M^0_5$}};
\fill (-1,1) circle (2pt) node[below] {\scriptsize{1}};
\fill (0,1) circle (2pt) node[below] {\scriptsize{3}};
\fill (1,1) circle (2pt) node[below] {\scriptsize{4}};
\draw[->, thin] (1,1) .. controls +(up:0.5cm)  .. (0,1);
\fill (2,1) circle (2pt) node[below] {\scriptsize{2}};
\draw[->, thin] (2,1) .. controls +(up:0.5cm)  .. (-1,1);
\fill (3,1) circle (2pt) node[below] {\scriptsize{3}};
\fill (4,1) circle (2pt) node[below] {\scriptsize{1}};
\fill (5,1) node[right] {\small{$\mathbf M^0_6$}};
\draw[->] (-3.5,0) -- (5.5,0);
\fill (-1,0) node[]{$|$} node[below]{$C_0=0$};
\fill (3,0) node[]{$|$} node[below]{$C_1$};
\end{tikzpicture}
\caption[smallcaption]{Construction of the increasing sequence 
$\suite{\mathbf M^0_n}$ and the constructions points, for  the compatibility graph $G$ of 
Figure \ref{fig:example1} and the sample $134231$. 
}
\label{fig:increase}
\end{center}
\end{figure} 

In particular, we can easily construct from the increasing sequence $\suite{\mathbf M^{0}_{n}}$, the limiting object 
$\mathbf M^{0}_{\infty}$ 
as the infinite labelled random graph obtained when letting $n$ go to infinity in the above. 
We then say that $\mathbf M^{0}_{\infty}$ is perfect if all its nodes are of degree 1. 
We have the following result, 

\begin{proposition}
If the graph $G$ is non-bipartite and condition (\ref{eq:Ncond}) is satisfied, then the infinite labelled random graph 
$\mathbf M^{0}_{\infty}$ is a.s. perfect. 
\end{proposition}

\begin{proof}
From Theorem \ref{thm:main}, the chain $\suite{W^{\{\emptyset\}}_n}$ 
is positive recurrent, so the family of integers 
\[\mathscr C:=\left\{n \in \N\,:\, W^{\{\emptyset\}}_n = \emptyset\right\}\]
is a.s. infinite. 
The elements of the latter set, denoted $0:=C_0 < C_1 < C_2 < ... $, are called {\em construction points} of the model. 
Then, we readily get \[\mathbf M^{0}_{\infty} = \lim_{k\to \infty} \mathbf M^0_{C_k}=\lim_{k\to\infty} \bigcup_{j=1}^{k} \mathbf M^{C_{j-1}}_{C_j} = \bigcup_{j=1}^{\infty} \mathbf M^{C_{j-1}}_{C_j},\]
which is a union of finite perfect {\sc fcfm}-matchings. It follows that all the nodes are of degree 1. 
\end{proof}

\subsection{Perfect {\sc fcfm}-matchings in reverse time} 
Perfect {\sc fcfm}-matchings have an interesting time-reversibility property. 
First, observe that we can complete the ``exchange" mechanism introduced in Section \ref{subsec:Markov}, using construction points as follows: 
Consider two construction points $C_m<C_{n} \in\mathscr C$. 
Then, in the perfect finite {\sc fcfm}-matching $\mathbf M^{C_m}_{C_{n}}$, replace the indexes of all nodes (i.e. the classes of the corresponding items) 
by the copies of the classes of their matches. 
Last, keep all edges of $\mathbf M^{C_m}_{C_{n}}$ unchanged. 
This procedure is illustrated in Figure \ref{Fig:exchange}. 
\begin{figure}[h!]
\begin{center}
\begin{tikzpicture}[scale=0.8]
%
\fill (0,1) circle (2pt) node[below] {\scriptsize{1}};
\draw[-, thin] (0,1) .. controls +(up:0.5cm)  .. (3,1);
\fill (1,1) circle (2pt) node[below] {\scriptsize{3}};
\draw[-, thin] (1,1) .. controls +(up:0.5cm)  .. (2,1);
\fill (2,1) circle (2pt) node[below] {\scriptsize{4}};
\fill (3,1) circle (2pt) node[below] {\scriptsize{$2$}};
\fill (4,1) circle (2pt) node[below] {\scriptsize{3}};
\draw[-, thin] (4,1) .. controls +(up:0.5cm)  .. (7,1);
\fill (5,1) circle (2pt) node[below] {\scriptsize{1}};
\draw[-, thin] (5,1) .. controls +(up:0.5cm)  .. (8,1);
\fill (6,1) circle (2pt) node[below] {\scriptsize{3}};
\draw[-, thin] (6,1) .. controls +(up:0.5cm)  .. (10,1);
\fill (7,1) circle (2pt) node[below] {\scriptsize{2}};
\fill (8,1) circle (2pt) node[below] {\scriptsize{2}};
\fill (9,1) circle (2pt) node[below] {\scriptsize{1}};
\draw[-, thin] (9,1) .. controls +(up:0.5cm)  .. (11,1);
\fill (10,1) circle (2pt) node[below] {\scriptsize{4}};
\fill (11,1) circle (2pt) node[below] {\scriptsize{2}};
\fill (13,1) node[right] {\small{$\mathbf M^{C_{n-2}}_{C_{n}}$}};
\fill (0,-0.5) circle (2pt) node[below] {\scriptsize{$\bar 2$}};
\draw[-, thin] (0,-0.5) .. controls +(up:0.5cm)  .. (3,-0.5);
\fill (1,-0.5) circle (2pt) node[below] {\scriptsize{$\bar 4$}};
\draw[-, thin] (1,-0.5) .. controls +(up:0.5cm)  .. (2,-0.5);
\fill (2,-0.5) circle (2pt) node[below] {\scriptsize{$\bar 3$}};
\fill (3,-0.5) circle (2pt) node[below] {\scriptsize{$\bar 1$}};
\fill (4,-0.5) circle (2pt) node[below] {\scriptsize{$\bar 2$}};
\draw[-, thin] (4,-0.5) .. controls +(up:0.5cm)  .. (7,-0.5);
\fill (5,-0.5) circle (2pt) node[below] {\scriptsize{$\bar 2$}};
\draw[-, thin] (5,-0.5) .. controls +(up:0.5cm)  .. (8,-0.5);
\fill (6,-0.5) circle (2pt) node[below] {\scriptsize{$\bar 4$}};
\draw[-, thin] (6,-0.5) .. controls +(up:0.5cm)  .. (10,-0.5);
\fill (7,-0.5) circle (2pt) node[below] {\scriptsize{$\bar 3$}};
\fill (8,-0.5) circle (2pt) node[below] {\scriptsize{$\bar 1$}};
\fill (9,-0.5) circle (2pt) node[below] {\scriptsize{$\bar 2$}};
\draw[-, thin] (9,-0.5) .. controls +(up:0.5cm)  .. (11,-0.5);
\fill (10,-0.5) circle (2pt) node[below] {\scriptsize{$\bar 3$}};
\fill (11,-0.5) circle (2pt) node[below] {\scriptsize{$\bar 1$}};
%
\fill (0,-2) circle (2pt) node[below] {\scriptsize{$\bar 1$}};
\draw[-, thin] (0,-2) .. controls +(up:0.5cm)  .. (2,-2);
\fill (1,-2) circle (2pt) node[below] {\scriptsize{$\bar 3$}};
\draw[-, thin] (1,-2) .. controls +(up:0.5cm)  .. (5,-2);
\fill (2,-2) circle (2pt) node[below] {\scriptsize{$\bar 2$}};
\fill (3,-2) circle (2pt) node[below] {\scriptsize{$\bar 1$}};
\draw[-, thin] (3,-2) .. controls +(up:0.5cm)  .. (6,-2);
\fill (4,-2) circle (2pt) node[below] {\scriptsize{$\bar 3$}};
\draw[-, thin] (4,-2) .. controls +(up:0.5cm)  .. (7,-2);
\fill (5,-2) circle (2pt) node[below] {\scriptsize{$\bar 4$}};
\fill (6,-2) circle (2pt) node[below] {\scriptsize{$\bar 2$}};
\fill (7,-2) circle (2pt) node[below] {\scriptsize{$\bar 2$}};
\fill (8,-2) circle (2pt) node[below] {\scriptsize{$\bar 1$}};
\fill (9,-2) circle (2pt) node[below] {\scriptsize{$\bar 3$}};
\draw[-, thin] (9,-2) .. controls +(up:0.5cm)  .. (10,-2);
\fill (10,-2) circle (2pt) node[below] {\scriptsize{$\bar 4$}};
\fill (11,-2) circle (2pt) node[below] {\scriptsize{$\bar 2$}};
\draw[-, thin] (8,-2) .. controls +(up:0.5cm)  .. (11,-2);
\fill (13,-2) node[right] {\small{$\cv{\td{\mbox{$\mathbf{M}$}^{C_{n}}_{C_{n-2}}}}$}};
\end{tikzpicture}
\caption[smallcaption]{Top: A perfect {\sc fcfm}-matching on two construction points, say $\mathbf M^{C_{n-2}}_{C_{n}}$.\\ 
Middle: The matching $\mathbf M^{C_{n-2}}_{C_{n}}$ after completing the exchanges on the two perfectly matched blocks. \\
Bottom: The resulting matching $\cv{\td{\mbox{$\mathbf{M}$}^{C_{n}}_{C_{n-2}}}}$ after time reversal is a {\sc fcfm}-matching.} 
\label{Fig:exchange}
\end{center}
\end{figure}

It is then immediate to observe that the resultant indexed random graph, which we denote $\cv{\td{\mbox{$\mathbf{M}$}^{C_{n}}_{C_m}}}$, 
is a perfect matching in the sense that all its nodes are a.s. of degree 1. We call it the {\em reversed} perfect matching between $C_n$ and $C_m$.  
We have the following result, 
\begin{proposition}
\label{pro:exchange}
For all $n>0$, the reversed perfect matching $\cv{\td{\mbox{$\mathbf{M}$}^{C_{n}}_0}}=\cv{\td{\mbox{$\mathbf{M}$}^{C_{n}}_{C_0}}}$ obtained in the above procedure  
is a perfect {\sc fcfm}-matching. 
\end{proposition}

\begin{proof}
Let the four nodes $i$, $j$, $k$ and $\ell$ be such that in $G$, $i \v k$, $j \v k$ and $i \v \ell$, 
and suppose that, for $m\le n$, after the exchange procedure over the matching $\mathbf M^{C_{m-1}}_{C_m}$, four copies $\td i$, $\td j$, $\td k$ and $\td{\ell}$ are read in that order, in reverse time, i.e. from right to left.  Let us also assume that the {\sc fcfm} rule in reverse time is violated on this quadruple: then 
the $\td k$ item is matched with the $\td j$ item while the $\td i$ item is still unmatched, and then the latter item is matched with the $\td{\ell}$ item. This occurs if and only if, in direct time, the four items of 
classes $i$, $j$, $k$ and $\ell$ arrive in that order, and the $k$ item choses the $j$ item over the $i$ item for its match, and then the unmatched $i$ item is matched with the $\ell$ item. 
This violates in turn the {\sc fcfm} policy, according to which the $k$ item should have been matched with the $i$ item instead of the $j$ item. Thus, 
the reversed perfect matching $\cv{\td{\mbox{$\mathbf{M}$}^{C_{m}}_{C_{m-1}}}}$ is a perfect {\sc fcfm}-matching. We finally observe that, by the very definition of the construction points, 
\begin{equation*}
\cv{\td{\mbox{$\mathbf{M}$}_0 ^ {C_{n}}}} = \bigcup_{m=1}^n \cv{\td{\mbox{$\mathbf{M}$}^{C_{m}}_{C_{m-1}}}}.
\end{equation*} 
It completes the proof. 
\end{proof} 

\medskip
\noindent 
\textbf{Acknowledgements.} The authors would like to warmly thank Jocelyn Begeot, at Universit\'e de Lorraine, for his comments and his careful reading of this work.

\providecommand{\bysame}{\leavevmode\hbox to3em{\hrulefill}\thinspace}
\providecommand{\MR}{\relax\ifhmode\unskip\space\fi MR }
\providecommand{\MRhref}[2]{%
  \href{http://www.ams.org/mathscinet-getitem?mr=#1}{#2}
}
\providecommand{\href}[2]{#2}

\end{document}